\documentclass{amsart}

\usepackage[latin1]{inputenc}
\usepackage{harvard}
\usepackage{graphicx}
\usepackage{psfrag}
\usepackage{amsthm,amsfonts,amssymb,amscd,amsmath}

\newtheorem{theorem}{Theorem}[section]

\newtheorem{lemma}[theorem]{Lemma}

\newtheorem{remark}[theorem]{Remark}

\newtheorem{definition}[theorem]{Definition}
\newtheorem{conjecture}[theorem]{Conjecture}

\newcommand{\field}[1]{\mathbb{#1}}

\def\C{\mathbb{C}}

\def\eps{\epsilon}

\def\fT{\field{T}}
\def\fR{\field{R}}

\def\cA{\mathcal{A}}
\def\cB{\mathcal{B}}
\def\cC{\mathcal{C}}
\def\cD{\mathcal{D}}

\def\cF{\mathcal{F}}

\def\cI{\mathcal{I}}

\def\cR{\mathcal{R}}

\def\cW{\mathcal{W}}

\def\cO{\mathcal{O}}

\date{2011-07-12}

\begin{document}
\title[Numerical Estimates]{A numerical  study of infinitely renormalizable area-preserving maps}
\author{Denis Gaidashev}
\address{
Department of Mathematics, Uppsala University, Uppsala, Sweden\\
{\tt gaidash@math.uu.se}}

\author{Tomas Johnson}
\address{
Department of Mathematics, Cornell University, Ithaca, NY 14853, USA\\
{\tt tomas.johnson@cornell.edu}
}

\begin{abstract}
It has been shown in \cite{GJ2} and \cite{GJM} that  infinitely renormalizable area-preserving maps admit  invariant Cantor sets with a maximal Lyapunov exponent equal to zero. Furthermore, the dynamics on these Cantor sets for any two infinitely renormalizable  maps is conjugated by  a transformation that extends to a differentiable function whose derivative is H\"older continuous of exponent $\alpha>0$.

In this paper we investigate numerically the specific value of $\alpha$. We also present numerical evidence that the normalized derivative cocycle with the base dynamics in  the Cantor set is ergodic.  Finally, we compute renormalization eigenvalues to a high accuracy to support a conjecture that the renormalization spectrum is real.
\end{abstract}

\maketitle

\setcounter{page}{1}

%\tableofcontents
\section{Introduction}
Numerical investigations of universality properties in dynamics have been historically important, and, in general, motivated the renormalization approach to universality. Such studies were behind the pioneering discovery of the Feigenbaum-Coullet-Tresser period doubling universality in unimodal maps \cite{Fei1}, \cite{Fei2}, \cite{TC}. % Since then universality ---  independence of the quantifiers of the geometry of orbits and bifurcation cascades in families of maps of the choice of a particular family --- has been demonstrated to be a rather generic phenomenon in dynamics \cite{Lyu}. 

Universality problems are typically approached via {\it renormalization}. In a renormalization setting one introduces  a {\it renormalization} operator on a functional space, and demonstrates that this operator has a {\it hyperbolic fixed point}. Period-doubling renormalization for dissipative two-dimensional maps has been extensively studied in \cite{CEK1,dCLM,LM}.  Compared to the Feigenbaum-Collet-Tresser one-dimensional  renormalization, the new striking feature of the two dimensional renormalization for highly dissipative maps is that {\it the restriction of the dynamics to their Cantor attractors  is not rigid}. Indeed, if the average Jacobians of $F$ and $G$ are different, for example, $b_F<b_G$, then the conjugacy $F \arrowvert_{\cC_F}  \, { \approx \atop { h} }    \, G \arrowvert_{\cC_G}$ is not smooth, rather it is a  H\"older continuous function with a definite upper bound on the H\"older exponent:
$\alpha \le {1 \over 2} \left(1+{\log b_G \over \log b_F }\right) <1.$ It turns out that the period-doubling renormalization for area-preserving maps is very different from the dissipative case.

A universal period-doubling cascade in families of area-preserving maps was observed by several authors in the early 80's \cite{DP,Hel,BCGG,Bou,CEK2,EKW1}.  The existence of a hyperbolic fixed point for the period-doubling renormalization operator 
$$R_{EKW}[F]=\Lambda^{-1}_F \circ F \circ F  \circ \Lambda_F,$$
where $\Lambda_F(x,u)=(\lambda_F x, \mu_F u)$ is an $F$-dependent {\it linear} change of coordinates, has been proved with computer-assistance in \cite{EKW2}, and in \cite{GJM}. The dynamics of the area-preserving fixed point map has been studied in \cite{GJ1,GJ2,J1,GJM}

We have proved in \cite{GJ2}  that {\it infinitely renormalizable} maps in a neighborhood of the fixed point of \cite{EKW2} admit  a ``stable'' Cantor set, that is the set on which the maximal Lyapunov exponent is zero. We have also shown in the same publication  that the conjugacy of stable dynamics is at least bi-Lipschitz on a submanifold of locally infinitely renormalizable maps  of a finite codimension. This has been further improved in \cite{GJM}, where it is shown that the conjugacy of dynamics on the Cantor set in fact  extends to a differentiable map whose  derivative is H\"older continuous with an exponent $\alpha>0$. This seems to be in stark contrast with the case of dissipative maps.

In this paper we continue the study of the dynamics on the stable Cantor set, and we provide numerical evidence in support of the following two conjectures that will be described in detail in the paper: 
{\it
\bigskip

\noindent i)  the normalized derivative cocyle on the Cantor set is ergodic with respect to an absolutely continuous invariant measure on the  circle;

\bigskip

\noindent ii) the H\"older exponent $\alpha$ satisfies $\alpha>0.070347$  for all infinitely renormalizable maps in a neighborhood of the fixed point.}

\bigskip

We also provide numerical evidence for the following conjecture

\bigskip
{\it
\noindent iii) the area-preserving period doubling operator has a real spectrum.}

\bigskip

The last conjecture can be compared with the corresponding conjecture in the one-dimensional case \cite{CCR}. In particular, we compute approximations to the leading $100$ eigenvalues of the operator, and they appear to be real, up to numerical accuracy.

\section{Renormalization for area-preserving reversible twist maps} 

An ``area-preserving map'' will mean an exact symplectic diffeomorphism of a subset of ${\fR}^2$ onto its image.

A reversible area-preserving map that satisfies a twist condition
$$\partial_u \left( \pi_x F(x,u) \right) \ne 0$$
everywhere in its domain of definition can be uniquely specified by a generating function $s$:

\begin{equation}\label{sdef}
\left({x  \atop  -s(y,x)} \right)  {{ \mbox{{\small \it  F}} \atop \mapsto} \atop \phantom{\mbox{\tiny .}}} \left({y \atop s(x,y) }\right).
\end{equation}

% Furthermore, the $2^k$ periodic orbits scale asymptotically with two scaling parameters
% \begin{equation}
% \lambda=-0.249 \ldots,\quad \mu=0.061 \ldots
% \end{equation}

The universality phenomenon in area-preserving maps can be explained rigorously if one shows that the following {\it renormalization} operator
\begin{equation}\label{Ren}
R_{EKW}[F]=\Lambda^{-1}_F \circ F \circ F \circ \Lambda_F,
\end{equation}
where $\Lambda_F$ is some $F$-dependent coordinate transformation, has a fixed point, and the derivative of this operator is hyperbolic at this fixed point. 
% $$
%  \left({x'  \atop  -s(Z,x')} \right) {{ \mbox{{\small \it  F}} \atop \mapsto} \atop \phantom{\mbox{\tiny .}}} \left({Z \atop s(x',Z)} \right) = \left({Z  \atop  -s(y',Z)} \right) {{ \mbox{{\small \it  F}} \atop \mapsto} \atop \phantom{\mbox{\tiny .}}} \left({y'\atop  s(Z,y')} \right).
% $$
According to \cite{CEK2}  $\Lambda_F$ can be chosen to be a linear diagonal transformation:  

$$\Lambda_F(x,u)=(\lambda x, \mu u).$$
The induced operator on generating functions will be denoted by $\cR_{EKW}$.
%We, therefore, set  $(x',y')=(\lambda x,  \lambda y)$, $Z(\lambda x, \lambda y)= z(x,y)$ to obtain:

% \begin{equation}\label{doubling}
% \left(\!{x  \atop  -{ 1 \over \mu } s(z,\lambda x)} \!\right) \!{{ \mbox{{\small $\Lambda_F$}} \atop \mapsto} \atop \phantom{\mbox{\tiny .}}} \!\left(\!{\lambda x  \atop  -s(z,\lambda x)} \!\right) \!{{ \mbox{{\small \it  F $ \circ$ F}} \atop \mapsto} \atop \phantom{\mbox{\tiny .}}}\!\left(\!{\lambda y \atop s(z,\lambda y)}\! \right)   {{ \mbox{{\small \it  $\Lambda_F^{-1}$}} \atop \mapsto} \atop \phantom{\mbox{\tiny .}}} \left(\!{y \atop {1 \over \mu } s(z,\lambda y) }\!\right),
% \end{equation}
% where $z(x,y)$ solves
% \begin{equation}\label{midpoint}
% s(\lambda x, z(x,y))+s(\lambda y, z(x,y))=0.
% \end{equation}

% If the solution of $(\ref{midpoint})$ is unique, then $z(x,y)=z(y,x)$, and it follows from $(\ref{doubling})$ that the generating function of the renormalized $F$ is given by 
% \begin{equation}
% \tilde{s}(x,y)=\mu^{-1} s(z(x,y),\lambda y).
% \end{equation}

%One can fix a set of normalization conditions for $\tilde{s}$ and $z$ which serve to determine scalings $\lambda$ and $\mu$ as functions of $s$. 
%For example, the normalization $s(1,0)=0$ is reproduced for $\tilde{s}$ as long as $z(1,0)=z(0,1)=1.$ In particular, this implies that 
% $$s(Z(\lambda,0),0)=0,$$
% which serves as an equation for $\lambda$.  Furthermore, the condition $\partial_1 s(1,0)=1$ is reproduced as long as $\mu=\partial_1 z (1,0).$

We will now summarize the definitions and theorems for the renormalization operator acting on generating functions originally due to the authors of \cite{EKW1} and \cite{EKW2}. Part $(i)-(iii)$ of Theorem \ref{Main_Theorem_1} were proved in \cite{EKW2}, Theorem \ref{Main_Theorem_2} was proved in \cite{GJ2}, and Theorem \ref{Main_Theorem_1} part $(iv)$ and Theorem \ref{Main_Theorem_3} were proved in \cite{GJM}.

Consider  the dyadic group,
\begin{equation}\label{dyadic_group}
\{0,1\}^\infty=\underleftarrow{\lim}\{0,1\}^n,
\end{equation}
where $\underleftarrow{\lim}$ stands for the inverse limit. An element $w$ of the dyadic group can be represented as a formal power series $w \rightarrow \sum_{k=0}^\infty w_{k+1} 2^k$. The {\it odometer}, or the {\it adding machine}, $p: \{0,1\}^\infty \rightarrow  \{0,1\}^\infty$ is the operation of adding $1$ in this group.

We will require a couple of more definitions for a precise statement of our results from \cite{GJM}.

\begin{definition}\label{B_space}
The Banach space of functions  $s(x,y)=\sum_{i,j=0}^{\infty}c_{i j} (x-\beta)^i (y-\beta)^j$, analytic on a bi-disk
$$\cD_\rho(\beta)=\{(x,y) \in \field{C}^2: |x-\beta|<\rho, |y-\beta|<\rho\},$$
for which the norm
$$\|s\|_\rho=\sum_{i,j=0}^{\infty}|c_{i j}|\rho^{i+j}$$
is finite, will be referred to as $\cA^\beta(\rho)$.

$\cA_s^\beta(\rho)$ will denote its symmetric subspace $\{s\in\cA^\beta(\rho) : s_1(x,y)=s_1(y, x)\}$.

We will use the simplified notation $\cA(\rho)$  and  $\cA_s(\rho)$ for  $\cA^0(\rho)$ and $\cA_s^0(\rho)$, respectively. 
\end{definition}

\begin{definition}\label{Fdefr}
The set of reversible twist maps $F$ of the form (\ref{sdef})  with $s \in \cB_\varrho(\tilde{s}) \subset \cA_s^\beta(\rho)$  will be referred to as $\cF^{\beta,\rho}_\varrho(\tilde{s})$:
\begin{equation}\label{F-set}
\cF^{\beta,\rho}_\varrho(\tilde{s})=\left\{ F: (x,-s(y,x)) \mapsto (y,s(x,y))| \quad s \in \cB_\varrho(\tilde{s}) \subset \cA_s^{\beta}(\rho)\right\}.
\end{equation}

We will also use the notation 
$$\cF^{\rho}_\varrho(\tilde{s}) \equiv \cF^{0,\rho}_\varrho(\tilde{s}).$$
\end{definition}

We are now ready to state our main theorems from \cite{GJM}. The original formulation of Main Theorem \ref{Main_Theorem_1} in \cite{EKW2} was in the space $\cA^{0.5}_s(1.6)$. The statement below was proved in \cite{GJM}.

\begin{theorem}\label{Main_Theorem_1} (Existence and Spectral properties)
There exists a polynomial $s_0: \field{C}^2 \mapsto \field{C}$, such that 
\bigskip 
\begin{itemize}
\item[i)] The operator $\cR_{EKW}$ is well-defined, analytic and compact in $\cB_{\varrho_0}(s_0) \subset \cA_s(\rho)$, with
$$\rho=1.75, \quad \varrho_0=5.79833984375 \times 10^{-4}.$$

\bigskip 

\item[ii)] There exists a function $s^* \in \cB_r(s_0) \subset \cA_s(\rho)$ with
$$r=1.1 \times 10^{-10},$$
such that
$$\cR_{EKW}[s^*]=s^*.$$ 

\bigskip 

\item[iii)] The linear operator  $D \cR_{EKW}[s^{*}]$  has two eigenvalues outside of the unit circle:
$$ 8.72021484375 \le \delta_1 \le 8.72216796875, \quad  \delta_2={1 \over \lambda_*},$$
where 
$$  -0.248875313689    \le \lambda_* \le  -0.248886108398438.$$
\bigskip

\item[iv)] The complement of these two eigenvalues in the spectrum is compactly contained in the unit disk. The largest eigenvalue in the unit disk is equal to $\lambda_*$, while
%\begin{eqnarray}
%\nonumber &&{\rm spec}(D {\cR}_{EKW}[s^*]) \setminus \{\delta_1,\delta_2 \} \subset \{z \in \C: |z| \le |\lambda_*|\},\\
$$
{\rm spec}(D \cR_{EKW}[s^*]) \setminus \{\delta_1,\delta_2, \lambda_*\} \subset \{z \in \C: |z| \le 0.1258544921875 \equiv \nu\}.
$$
%\end{eqnarray}

\end{itemize}
\end{theorem}

The Main Theorem $\ref{Main_Theorem_1}$ imply that there exist codimension $2$ local stable manifolds $\cW_{\cR_{EKW}}(s^*) \subset \cA_s(1.75)$ of the operator $\cR_{EKW}$. %and  $\cW_{\cR_{EKW}}^{0.5}(s^{EKW}) \subset \cA_s^{0.5}(1.6)$ of the operator $\cR_{EKW}$. 

Compactness of the operator  $\cR_{EKW}$ in  neighborhood of $s^*$ implies that there exists a strong ``submanifold'' 
$$\cW_{\cR_{EKW}}^s(s^*) \subset \cW_{\cR_{EKW}}(s^*),$$
of codimension $1$ in $\cW_{\cR_{EKW}}(s^*)$, such that the contraction rate in $\cW_{\cR_{EKW}}^s(s^*)$ is bounded from above by the number $\nu$ from Part $iv)$ from $\ref{Main_Theorem_1}$:
$$\|\cR_{EKW}^n[s]-\cR_{EKW}^n[\tilde{s}]\|_\rho=O(\nu^n)$$
for any two $s$ and $\tilde{s}$ in $\cW_{\cR_{EKW}}^s(s^*)$.

\begin{definition}\label{Wdeff}
The set of reversible twist maps of the form (\ref{sdef}) such that $s \in \cW_{\cR_{EKW}}(s^*) \subset \cA_s(1.75)$ will be denoted $W_{EKW}$, and referred to as  {\bf{infinitely renormalizable}} maps.

The set of reversible twist maps of the form (\ref{sdef}) such that $s \in \cW^s_{\cR_{EKW}}(s^*) \subset \cA_s(1.75)$ will be denoted $W_{EKW}^s$.

%The set of reversible twist maps of the form (\ref{sdef}) such that $s \in \cW_{\cR_{EKW}}^{0.5}(s^{EKW}) \subset \cA_s^{0.5}(1.6)$ will be denoted $W_{EKW}^{0.5}$. 
\end{definition}

%Recall the Definition $\ref{Fdefr}$.

\begin{definition}\label{Wdefr}
Set,
\begin{eqnarray}
%\nonumber \cW_\varrho(\tilde{s}) &\equiv&  \cW \cap \cF^{1.75}_\varrho(\tilde{s}), \\
\nonumber W_\varrho(\tilde{s}) &\equiv&  W_{EKW} \cap \cF^{1.75}_\varrho(\tilde{s}),\\
\nonumber W^s_\varrho(\tilde{s}) &\equiv&  W_{EKW}^s \cap \cF^{1.75}_\varrho(\tilde{s}).
%\nonumber W^{0.5}_\varrho(\tilde{s}) &\equiv&  W_{EKW}^{0.5} \cap \cF^{0.5,1.6}_\varrho(\tilde{s}).
\end{eqnarray}

\end{definition}

Naturally, these sets are  invariant under renormalization if $\varrho$ is sufficiently small.

The eigenvector of $\lambda_*$ is denoted by $\psi_{s^*}^{EKW}=\psi_{s*}+\tilde \psi$ and $h_{\lambda_*}^{EKW}$, for $D \cR_{EKW}$ and $D R_{EKW}$, respectively ( see \cite{GJM} for details), where $\psi_{s*}$ is the eigenvector with eigenvalue $\lambda_*$ for the operator $D\cR_0$ induced by $R_0 [F] = \Lambda_*^{-1} \circ F\circ F \circ \Lambda_*$. $\tilde \psi$ is an eigenvector of $D\cR_0$ and $D\cR_{EKW}$ with eigenvalue $1$ and $0$, respectively. The $\lambda_*$ direction is constructed by a coordinate change of the form $S_\epsilon^{-1}\circ F \circ S_\epsilon$, where
\begin{eqnarray}
\label{coord} S_\epsilon(x,u)&=&\left(x+\epsilon x^2,{u \over 1+2 \epsilon x}\right),\\
%\label{coord} &=&(x,u) +\epsilon h_{\lambda_*}^{EKW}+O(\epsilon^2),\\
\label{inv_coord} S_\epsilon^{-1}(y,v)&=&\left({\sqrt{1+4\epsilon y}-1 \over 2 \epsilon}, v \sqrt{1+4 \epsilon y}  \right).
\end{eqnarray}
The map $S_\epsilon^{-1}\circ F \circ S_\epsilon$  is generated by:
$$\hat{s}(x,y)=s(x+\epsilon x^2, y+ \epsilon y^2) (1+2 \epsilon y).$$

In addition, the following more general fact about the spectrum of the operator $\cR_0[s^*]$ holds (see  \cite{GJM}). Denote by $\cO_2(\cD)$ the Banach space of maps $F: \cD \mapsto \field{C}^2$, analytic on an open simply connected set $\cD \subset \field{C}^2$, continuous on $\partial \cD$, equipped with a finite max supremum norm $\| \cdot \|_\cD$:
$$\| F \|_{\cD}=\max\left\{\sup_{(x,u) \in \cD}|F_1(x,u)|, \sup_{(x,u) \in \cD}|F_2(x,u)| \right\}.$$

\begin{lemma}\label{powers}
There exists a domain $\cD$, such that the operator $R_0$ is differentiable in a neighborhood of $F^*$ in $\cO_2(\cD)$, and ${\rm spec}(D R_0[F^*])$ contains eigenvalues $\mu_*^j \lambda^i_*$, $i \ge 0$, $j \ge 0$, of multiplicity at least $2$, and $\mu_*^j \lambda_*^{-1}$, $j \ge 0$, and $\lambda_*^i \mu^{-1}_*$, $i \ge 0$, of multiplicity at least $1$.
\end{lemma} 

\begin{remark}
Note, that not all of the maps in $\cO_2(\cD)$ are generated by generating functions, neither do all of them admit the symmetry $T\circ  F \circ T=F^{-1}$, therefore not all of the products of powers of $\mu_*$ and $\lambda_*$ mentioned in Lemma $\ref{powers}$ are eigenvalues  of $D R_{EKW}[F^*]$ in $T_{F^*} \cF^{1.75}_\varrho(s^*)$.
\end{remark}

%The operators $\cR$ and $R$ are constructed from $\cR_{EKW}$ and $R_{EKW}$ by removing the $\lambda_*$ direction.

The following two theorems describe the dynamics of the infinitely-renormalizable maps on their ``stable'' Cantor sets.

\bigskip

\begin{theorem}\label{Main_Theorem_2}(Stable Set)

There exists $\varrho>0$ such that any $F \in W_\varrho(s_0)$,  admits a ``stable'' Cantor set $\cC_F \subset \cD$ with the following properties.
 
\begin{itemize}
\item[i)] For all $x \in \cC_F$ the maximal Lyapunov exponent $\chi(x;F)$ exists, is $F$-invariant, is equal to zero:
$$\chi(x;F) =0,$$
and 
$$\lim_{i \rightarrow \pm \infty} {1 \over |i|} \log\left\{  { \Arrowvert D F^i(x) v \Arrowvert  \over \|v\| }\right\}=0,$$
uniformly for all $v \in \fR^2 \setminus \{0\}$ and $x \in \cC_F$.

\item[ii)] The Hausdorff dimension of $\cC_F$ satisfies
$${\rm dim}_H(\cC_F) \le 0.794921875.$$

\item[iii)] The restriction of the dynamics $F \arrowvert_{\cC_F}$ is topologically conjugate to the adding machine.
\end{itemize}
\end{theorem}

\begin{theorem}\label{Main_Theorem_3}(Rigidity)
Let $s^*$ and $\cC_F$ be as in Main Theorems  \ref{Main_Theorem_1} and \ref{Main_Theorem_2}. There exists $\varrho>0$, such that for all  $F$ and $\tilde{F}$  in $W_\varrho(s^*)$, 
$$F \arrowvert_{\cC_F}  \, { \approx \atop { h} }    \, \tilde{F} \arrowvert_{\cC_{\tilde{F}}},$$
where $h$ extends to a neighborhood of $\cC_F$ as a differentiable transformation, whose derivative $D h$ is H\"older continuous with the H\"older exponent 
$$\alpha \ge 0.0129241943359375.$$
\end{theorem}

\section{H\"older continuity of the derivative}\label{SnumRig}
In this section we will give numerical evidence in support of the following conjecture. 
\begin{conjecture}\label{CR1}
For all $F, \tilde{F} \in W_{\varrho_0}(s^*)$, $\varrho_0 = 6 \times 10^{-7}$, 
$$F \arrowvert_{\cC_F}  \, { \approx \atop { h} }    \, \tilde{F} \arrowvert_{\cC_{\tilde{F}}},$$
where $h$ is  a differentiable transformation, whose derivative $D h$ is H\"older continuous with positive H\"older exponent larger than:
$$ \alpha \geq 0.070347$$ 
\end{conjecture}

We would like to emphasize that it is not the $C^{1+\alpha}$ property that we are after, but rather the specific value of $\alpha$. The fact that the conjugacy extends to a differentiable function whose derivative is H\"older has been proved in \cite{GJM}.

To investigate the H\"older exponent of the derivative of the conjugacies in $W_{\varrho_0}(s_0)$, we conduct an experiment, which will be described below. The convergence of $Dh_\omega$ depends on the spectral properties of the operator $DR_{EKW}$ in the codimension one submanifold $W^s_{\varrho_0}$, which is transversal to $h_{\lambda_*}^{EKW}$ at $F^*$. Note that all our computations are done with the entire neighborhood ($F\in W_{\varrho_0}$) of $F^*$ on the stable manifold. To ensure that we study two different functions, we study a $\tilde F\in W_{\varrho_0}$ transported along the one dimensional manifold tangent to the eigenvector $h_{\lambda_*}^{EKW}$ at its intersection with $W^s_{\varrho_0}$, which we denote by $I(\tilde F)$. In practice this means that we study conjugacies from an arbitrary $F\in W^s_{\varrho_0}$, to an arbitrary $\tilde F\in W^s_{\varrho_0}$, transported some distance along $I(\tilde F)$. We study such functions using the generating function $\hat
  s_\eps$, defined by

\begin{eqnarray}
\nonumber \hat s_\eps(x,y)&:=&s^*((1-\eps)x+\eps (1-\eps)^2x^2,(1-\eps)y+\eps (1-\eps)^2y^2)(1+2\eps(1-\eps) y)/(1+\eps) \\
\nonumber &=&s^*(x,y)+\eps\psi_{s^*}+\eps\tilde\psi+O(\eps^2).
\end{eqnarray} 

The corresponding scalings are denoted by $\lambda_\eps$ and $\mu_\eps$, they are computed from the formulae
\begin{equation}
\hat s_\eps(\lambda_\eps,1)+\hat s_\eps(0,1) = 0
\end{equation}
and
\begin{eqnarray}
\mu_\eps=\lambda_\eps\partial_1 \hat s_\eps(Z[\hat s_\eps](x,y),y)_{|_{(\lambda_\eps,0)}},
\end{eqnarray}
where $\partial_1 Z[\hat s_\eps](x,y)$ is computed using the midpoint formula:
$$
\hat s_\eps(x,Z[\hat s_\eps](x,y))+\hat s_\eps(y,Z[\hat s_\eps](x,y)) = 0, 
$$
i.e.,
\begin{eqnarray*}
0 & = &\partial_1 \left(\hat s_\eps(x,Z[\hat s_\eps](x,y))+\hat s_\eps(y,Z[\hat s_\eps](x,y))\right) \\
& = & \partial_1 \hat s_\eps(x,Z[\hat s_\eps](x,y))+\partial_2 \hat s_\eps(x,Z[\hat s_\eps](x,y))Z_1[\hat s_\eps](x,y) \\ & & + \, \partial_2 \hat s_\eps(y,Z[\hat s_\eps](x,y))Z_1[\hat s_\eps](x,y),
\end{eqnarray*}
from which
\begin{equation}
\partial_1 Z[\hat s_\eps](\lambda_\eps,0) = \frac{-\partial_1 \hat s_\eps(\lambda_\eps,1)}
{\partial_2 \hat s_\eps(\lambda_\eps,1)+ \partial_2 \hat s_\eps(0,1)}
\end{equation}
follows, and hence
\begin{eqnarray}
\mu_\eps=\frac{-\lambda_\eps \partial_1 \hat s_\eps(1,0) \partial_1 \hat s_\eps(\lambda_\eps,1)}
{\partial_2 \hat s_\eps(\lambda_\eps,1)+ \partial_2 \hat s_\eps(0,1)}.
\end{eqnarray}

To approximate the renormalization operator acting on $\hat s_\eps$, we note that 
\begin{eqnarray*}
\cR_{EKW}[\hat s_\eps]& = &\cR_{EKW}[s^*]+\eps D\cR_{EKW}[s^*]\psi_{s^*}+\eps D\cR_{EKW}[s^*]\tilde \psi+O(\eps^2)\\
& = &\cR_{EKW}[s^*]+\eps D\cR_{EKW}[s^*]\psi_{s^*}+O(\eps^2)\\
& = & s^*+\eps D\cR_{0}[s^*]\psi_{s^*}+\eps\lambda_*\tilde\psi+O(\eps^2) \\
& = & s^*+\eps \lambda_*\psi_{s^*}+\eps\lambda_*\tilde\psi+O(\eps^2) \\
& = & \hat s_{\eps\lambda_*}+O(\eps^2).
\end{eqnarray*}
Thus,
\begin{eqnarray*}
\cR^n_{EKW}[\hat s_\eps]& = & \hat s_{\eps(\lambda_*)^n}+O(\eps^2)+\cdots+O(\eps^2(\lambda_*)^{2(n-1)}) \\
& = & \hat s_{\eps(\lambda_*)^n}+O(\eps^2/(1-(\lambda_*)^2)) \\
& = & \hat s_{\eps(\lambda_*)^n}+O(\eps^2).
\end{eqnarray*}
We denote by $F_\eps$ the map generated by a generating function $\hat s_\eps$. To simplify the notation we denote the corresponding presentation functions by 

\begin{eqnarray*}
\Psi^\eps_0 & := & \Lambda_{F_\eps} \\
\Psi^\eps_1 & := & F_\eps\circ\Lambda_{F_\eps} \\
\Psi^\eps_{\omega^n} & := & \Psi^\eps_{\omega_1}\circ \cdots \circ \Psi^{\eps(\lambda_*)^{(n-1)}}_{\omega_n}.\\
\end{eqnarray*}

According to \cite{GJM}, the maps
$$h^\eps_{\omega^n} := \Psi^\eps_{\omega^n}\circ \left(\Psi^*_{\omega^n}\right)^{-1}
$$
converge to a conjugacy of  $F^*\arrowvert_{\cC_*}$ and  $F_\eps \arrowvert_{\cC_{F_\eps}}$, which is, furthermore, extendable to a $C^{1+\alpha}$ map on a neighborhood of $\cC_*$.

The following Lemma about the existence of hyperbolic fixed points for maps in a small neighborhood of the renormalization fixed point  map $F^*$ is a restatement of a result from \cite{GJ1} in the setting of the functional space $\cA_s(1.75)$. The proof of the Lemma is computer-assisted (see \cite{GJ1}).

\begin{lemma}\label{hyp_point}
Every map $F \in \cF^{\rho}_\varrho(s_0)$, with $\varrho=6.0 \times 10^{-12}$ and $\rho=1.75$, possesses a hyperbolic fixed point $p^F \in \cD$, such that
\begin{itemize}
\item[1)] $\pi_x p^F \in (0.577606201171875,0.577629923820496)$, and $\pi_u p^F=0$, where $\pi_{x,u}$ are projections on the $x$ and $u$ coordinates;
\item[2)] $D F(p^F)$ has two negative eigenvalues. 
\begin{eqnarray}
\nonumber e_+^F &\in& (-2.0576171875,-2.057373046875),\\
\nonumber e_-^F &\in& (-0.486053466796875,-0.48602294921875),
\end{eqnarray} 
corresponding to the following two eigenvectors:
$$
{\bf s}^F=[1.0,-(0.77978515625,0.779815673828125)], \quad {\rm and} \quad {\bf u}^F=T({\bf s}^F).
$$
\end{itemize}
\end{lemma}

\bigskip

This Lemma implies the existence of hyperbolic $2^n$-th periodic orbits for maps in $\cW_\varrho(s_0)$.  Let $\cO_n(F)$ denote such $2^n$-th periodic orbit  of $F \in  \cW_\varrho(s_0)$, specifically:
$$\cO_n(F)=\bigcup_{i=0}^{2^n-1} F^{i}(\Psi^{F}_{0^n}(p^{F_n})),$$
where $p^{F_n}$ is the fixed point of $F_n \equiv R^n[F]  \in \cW_\varrho(s_0)$. We will also denote 
$$ p_{0^n}^F=\Psi^{F}_{0^n}(p^{F_n}),\quad p_{\omega}^F= F^{\sum_{i=1}^n \omega_i 2^{i-1}}(p_{0^n}^F).$$

Consider the $2^n$-periodic point $p_{\omega^n} \equiv p_{\omega^n}^{F^*} \in B^*_{\omega^n}$.
% where we have set
%$$
%p_{\omega^n}=p_{\omega^n}^{F^*} = \Psi^*_{\omega^n}(p^{F^*}),
%$$
%and where $p^{F^*}$ is the fixed point of $F^{*}$ computed in \cite{GJ1}. 

We will estimate the H\"older exponent of the conjugacy $h^\eps$ by computing  perturbations of $p_{\omega^n}$ in $B^*_{\omega^n}$.

To this end, we estimate 
\begin{equation}
N_n(\epsilon, \omega^n, t) := 
\|Dh^\eps_{\omega^n}(p_{\omega^n}) - Dh^\eps_{\omega^n}(p_{\omega^n}+\delta(\omega^n)(\cos t,\sin t))\|,
\end{equation}
where, for $n\leq 15$, 
$$
\delta(\omega^n) := \theta^{2.042n} \|p_{\omega^n}\|,
$$
with $\theta=0.272$ turns out to be a good choice. The reason that we multiply with $\|p_{\omega^n}\|$ is that the size of different points on a periodic orbit can vary dramatically, and we want to make the perturbations on one renormalization level to be, relatively, of the same order, i.e., we perturb the same significant digit of each point on the periodic orbit. We remark that $\delta$ needs to be sufficiently small to guarantee that the perturbation is inside of $B^*_{\omega^n}$, and at the same time sufficiently large to be computable. In particular we want to be able to compute $N_n(\epsilon, \omega^n, t)$ with $n$ as large as possible. Our computations are done with IEEE extended precision (80 bits), for which the $\delta$ chosen above is a good compromise, and allows us to compute with $n\leq 15$.

We are interested in bounds on $\alpha$ that are uniform in $t$ and $\omega^n$. Therefore, we define
\begin{eqnarray}
M_n(\epsilon, \omega^n) & := & \max_{t\in [0,2\pi]} N_n(\epsilon, \omega^n, t) \\
\alpha_n(\epsilon) & := & \min_{\omega^n \in \{0,1\}^n}\frac{\log M_n(\epsilon, \omega^n)}{\log \delta(\omega^n)}.  
\end{eqnarray}
The motivation to the definition of $\alpha_n(\eps)$ is that we want to find an $\alpha$ such that 
$$
\|Dh^\eps_{\omega^n}(x) - Dh^\eps_{\omega^n}(y)\| \leq C \|x-y\|^\alpha,
$$
for all $x\neq y \in B^*_{\omega^n}$. 

As $n \rightarrow \infty$, $\delta(\omega^n)\rightarrow 0$ uniformly, since ${\rm diam}(B^*_{\omega^n})\leq C \theta^n$. Thus, we get the following dichotomy, if 
$$\alpha_n(\eps) \rightarrow \alpha(\eps) > 0,
$$
then $\alpha$ exists and is equal to $\alpha(\epsilon)$, and if there exists $m$, such that
$$
\alpha_n(\eps) < 0, \quad \textrm{for }n>m,
$$  
then $\alpha$ does not exist.%, i.e., $Dh^\eps_\omega$ is not H\"older on $\cI$, and hence not on $\cV_r(s_0)$. 

Obviously, we are only able to collect a finite amount of data, i.e., we are only able to compute $N_n$, for a finite discretization of $[0,2\pi]$, and only for a finite number of renormalization levels. The first discretization is not a major concern, since by discretizing until the results stabilize, we can be almost certain that our results are correct (up to some accuracy), this issue will be explained further in the section below. In the results given below we have discretized $[0,2\pi]$ with $512$ points. The computed values of $\alpha_n(\eps)$ agree to $5$ digits when the discretization is increased from $256$ to $512$ points. It is therefore reasonable to conclude that further discretization of the interval $[0,2\pi]$ would have a negligible effect on our results.

The second problem is much greater; ideally we should compute $\alpha_n(\eps)$ for large enough $n$'s so that $\alpha_n(\eps)$ converges. Unfortunately, this is not possible; on one hand due to the enormous complexity of the computation of $N_n$, which is $O(2^n)$, and on the other hand due to the finite accuracy of computers. Our choice of 
$\delta(\omega^n)=\theta^{2.042n}\|p_{\omega^n}\|$, 
means that $\delta(\omega^{16}) \approx 3.36\times10^{-19}$, and since the precision of the 80-bit extended precision is 19 decimals, this is already too small to be reliably computable. Note that for some choices of $\omega^n$,  $p_{\omega^n}^{F^*}$, might be macroscopic, regardless of $n$, preventing us from proceeding further. In principle we could use software implemented, higher precision libraries to compute for higher $n$'s, but, as we will show below, we would need to compute the first $30-50$ levels to get convergence, and $2^{30}$ is a prohibitively large number of combinations. Therefore, we stop at $n=15$.

In order to say something about the convergence of $\alpha_n(\eps)$ from a finite number of renormalization levels, i.e., from computing a finite number of $h^\eps_{\omega^n}$, we try to fit the results to a graph of a function which approaches a constant positive value. It turns out that for any choice of $\eps$ in the range, where our computation yields results, our measurements fit almost perfectly to a function $a_\eps(n)$ of the form 
$$
a_\eps(n):=a_\eps \exp{\left(k^1_\eps/n + k^2_\eps/n^2+k^3_\eps/n^3 + k^4_\eps/n^4\right)},
$$ 
clearly, $a_\eps(n)\rightarrow a_\eps$. We set 
$$
\alpha(\eps) := a_\eps.
$$

% \begin{figure}[h]
% \begin{center}
% \includegraphics[width=0.7 \textwidth]{alpha_comp_meas.eps}
% \caption{\textit{Tracking the measurements:} the blue curve is the measured value of $\alpha_n(\eps)$ as a function of the renormalization level. The green curve is the estimate of the form $a_\eps(n)=a_\eps \exp{\left(k/n+k^2_\eps/n^2+k^3_\eps/n^3 + k^4_\eps/n^4\right)}$. The results in this figure are from the computations with $\eps=1.593584796420859\times10^{-8}.$}\label{FalphaEst}
% \end{center}
%  \end{figure}

\begin{figure}[h]
\begin{center}
\includegraphics[width=0.49 \textwidth]{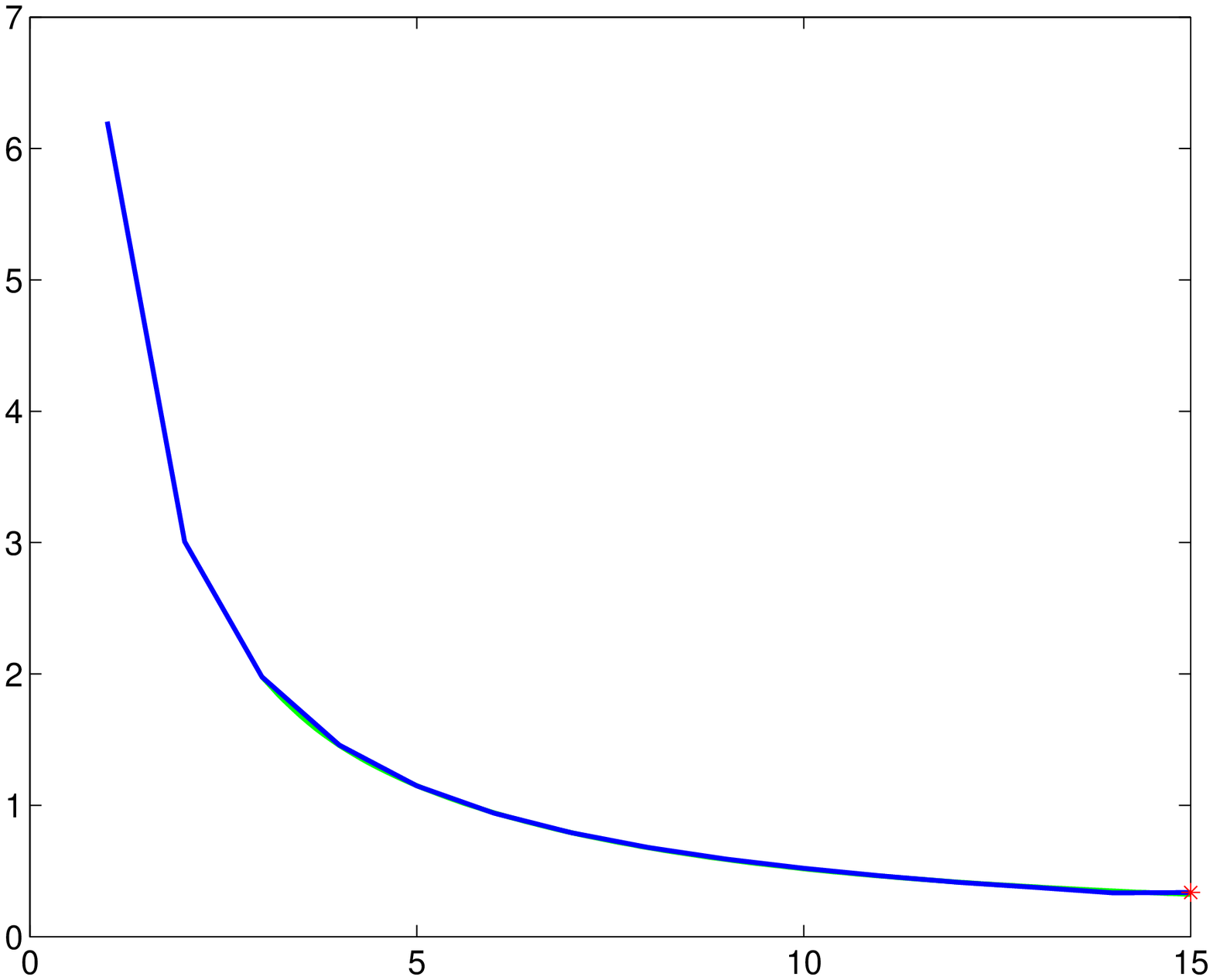}
\includegraphics[width=0.49 \textwidth]{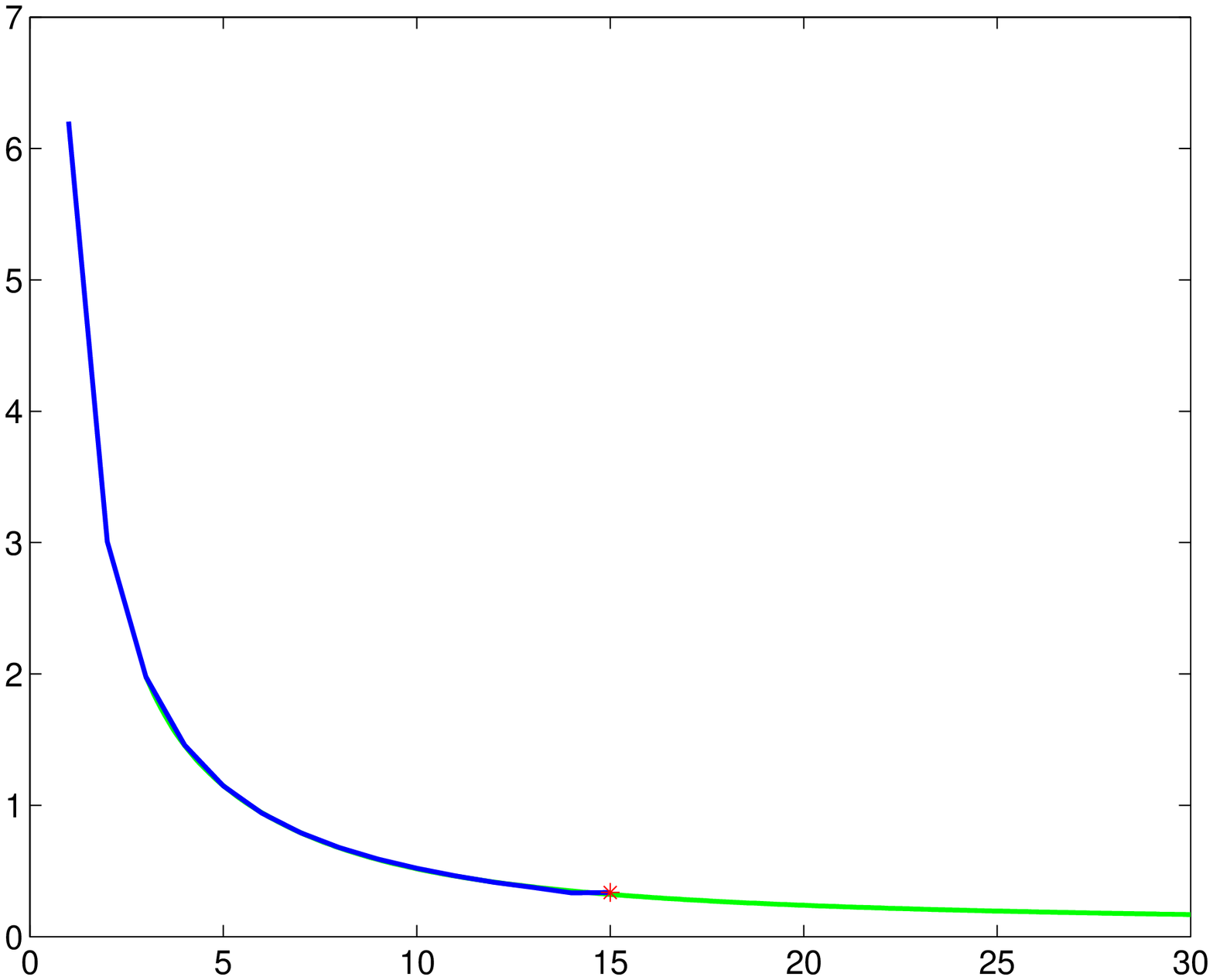}
\caption{(a)\textit{Tracking the measurements:} the blue curve is the measured value of $\alpha_n(\eps)$ as a function of the renormalization level. The green curve is the estimate of the form $a_\eps(n)=a_\eps \exp{\left(k/n+k^2_\eps/n^2+k^3_\eps/n^3 + k^4_\eps/n^4\right)}$. The results in this figure are from the computations with $\eps=1.593584796420859\times10^{-8}.$ (b) \textit{Extrapolation of the measurements:} to extrapolate the measurements, the measured data, $\alpha_n(\eps)$, is glued together with the graph of the computed estimate, $a_\eps(n)$. The blue curve is the measured value of $\alpha_n(\eps)$ as a function of the renormalization level, the red star is the last measurement. After the red star, the green curve is the graph of the function $a_\eps(n)=a_\eps \exp{\left(k/n+ k^2_\eps/n^2+k^3_\eps/n^3 + k^4_\eps/n^4\right)}$. The results in this figure are from the computations with $\eps=1.593584796420859\times10^{-8}.$}\label{FalphaExt}\label{FalphaEst}
\end{center}
 \end{figure}

This estimation of the measurements by the function $a_\eps(n)$ is illustrated in Figure \ref{FalphaEst}, where we have restricted to the subset $3\leq n\leq 15$, the reason being that we would like to avoid transient effects. These transient effects arise primarily from the fact that for small $n$, $\delta(\omega^n)$ is rather large compared to $M_n(\eps,\omega^n)$. 

The combined function constructed by gluing together the measurements $\alpha_n(\eps)$ and the graph of the function $a_\eps(n)$ at $n=15$ is illustrated in Figure \ref{FalphaExt}. From Figure \ref{FalphaExt} we see that the function does not flatten out completely (up to the resolution of the image) until $n=30$; and as already mentioned, computing $h_{\omega^{30}}^\eps$, is completely unrealistic. The computation at our final renormalization level, $n=15$, already takes about 4 hours for one choice of $\epsilon$.

We implement an algorithm based on the discussion above in a C++ program, \cite{JP}, and compute the estimated rigidity of $h_\omega^\eps$ for the largest possible range of $\epsilon$. This range turns out to be $[10^{-8}, 10^{-4}]$. The result of the computations together with the relative least squares error of the estimation
is given in Table \ref{TalphaOfEps}. 

\begin{table}[h]
\begin{center}
\begin{tabular}{r|l|r}
$\epsilon$ & $\alpha(\eps)$ & $\frac{1}{13}\sqrt{\sum_{n=3}^{15} | \frac{\alpha_n(\eps)-a_\eps(n)}{\alpha_n(\epsilon)} |^2}$\\ \hline
$10^{-4} $& $3.37\times10^{-3}$ & $0.40\%$\\
$10^{-5} $& $1.13\times10^{-2}$ & $0.25\%$\\
$10^{-6}$ & $2.21\times10^{-2}$ & $0.19\%$\\
$10^{-7}$ & $3.59\times10^{-2}$ & $0.13\%$\\
$1.593584796420859\times10^{-8}$ & $7.03\times10^{-2}$ & $0.50\%$
\end{tabular}
\end{center}
\caption{The estimated values of $\alpha(\epsilon)$ for various values of $\epsilon$, together with the least squares error of the estimating function $a_\eps(n)$.
}\label{TalphaOfEps}
\end{table}

To discuss how the rigidity of the conjugacies between $F_\eps$ and $F^*$ depends on the distance of $F_\eps$ from the fix point; we need to relate the size of $\eps$ with $\|\hat s_\eps - s^*\|$. This correspondence can be estimated using the fact that:
$$
\hat s_\eps-s^*=\eps\psi^{EKW}_{s^*}+O(\eps^2)
$$

Thus,
$$
\|\hat s_\eps-s^*\| = \eps\|\psi^{EKW}_{s^*}\| + O(\eps^2).
$$

By using the estimate 
$$
\|\psi^{EKW}_{s^*}\| \ge 37.6509616148184234,
$$
together with the estimate given in Table \ref{TalphaOfEps}; we conclude that $\varrho_0$ corresponds to an $\epsilon \le  1.593584796420859\times10^{-8}$ and that
$$
F^{\eps} \arrowvert_{\cC_{F^\eps}}  \, { \approx \atop { h} }    \, F^* \arrowvert_{\cC_*}, \quad \textrm{for  } \, F^\eps  \in \cI \subset  \cW_{\varrho_0}(s_0), \quad \textrm{with} \quad \alpha =  0.070347113466231. 
$$
This conclusion implies Conjecture \ref{CR1}.

%---------------------------------------------------------------------------------------------------------------------

\section{Twisting of line fields}\label{STLF}

Consider the stable and unstable invariant direction fields on the $2^n$-th periodic orbit $\cO_n(F)$. At every point $p_\omega^F$, $\omega \in \{0,1\}^n$  of $\cO_n(F)$, these directions are given by
\begin{eqnarray}
\label{vec_s} {\bf s}_\omega^F&=& D \Psi^F_\omega(p^{F_n}) {\bf s}^{F_n},\\
\label{vec_u} {\bf u}_\omega^F&=&D \Psi^F_\omega(p^{F_n}) {\bf u}^{F_n}.
\end{eqnarray}

The angles between these vectors and the positive real line will be denoted by $\alpha_\omega^F$ and  $\beta^F_\omega$, respectively.

\begin{figure}[t]
\begin{center}
\includegraphics[angle=-90,width=0.9 \textwidth]{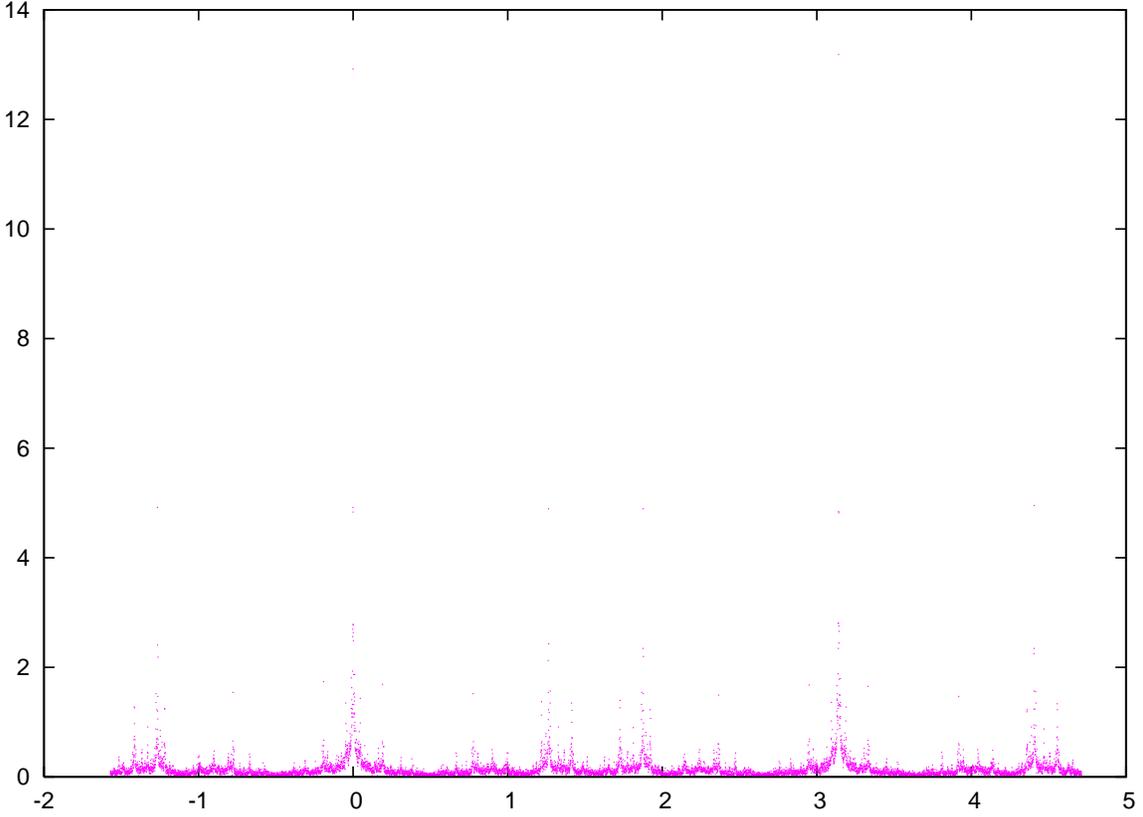}
\caption{
Distribution of angles in $\cO_{12}(F^*)$
}\label{st_angles_dist}
\end{center}
 \end{figure}

We have plotted the distribution of the angles $\alpha_\omega^F$ for $\omega \in \{0,1\}^{12}$ for  vectors  $(\ref{vec_s})$ in Fig. $\ref{st_angles_dist}$.

The figure demonstrates that the stable (and, via reversibility, the unstable) direction undergoes a sort of twisting: the directions do not concentrate, but are distributed with a certain measure over the whole circle.

One can heuristically explain the structure of this distribution. First, the appearance of two large peaks at angles $0$ and $\pi$ can be explained by  the fact that  all vectors ${\bf s}_\omega^F$ with $0$ in the beginning of their coding $\omega$ tend to be aligned horizontally: indeed, for every $\omega \in \{0,1\}^{n-k}$, 
$$\tan \alpha_{0^k \omega}^F= {\rm const} \ \left( {|\mu_*| \over \lambda_* }\right)^k \tan \alpha_\omega^{F_{k-1}}.$$
The smaller peaks (such as the two peaks between angles $1$ and $2$) correspond to the stable directions at points $p_\omega^F$ that are direct and inverse images under dynamics $F$ of points $p_\omega^F$ whose angles are in one of two central peaks.

Let a $v$ be a unit vector. We identify the set of unit vectors in $\fR^2$ with $\fT^1$, and write $v \in \fT^1$. Consider the normalized derivative cocycle $\cA_i: \cC_F \times \fT^1 \mapsto \fT^1$:
$$\cA_i(x) v \equiv {D F^i(x) v \over  \Arrowvert D F^i(x) v \Arrowvert}.$$

Fig. $\ref{st_angles_dist}$ lends evidence to the following conjecture:

%\begin{conjecture}
%There exists an absolutely continuous invariant measure for the cocycle $\cA_i$. %, with density as in \ref{st_angles_dist}
%\end{conjecture}

\begin{conjecture}
%$\phantom{aaa}\\$
%\begin{itemize}
%\item[1)] 
There exists an absolutely  continuous invariant measure $\nu$ on $\fT^1$, which is ergodic for the cocycle  $\cA_i(x)$, $x \in \cC_F$:
\begin{equation}
\label{ergodic}\lim_{n \rightarrow \infty} {1 \over n} \sum_{k=0}^{n-1} f \left( \cA_k(x) v \right) ={1 \over  \nu(\fT^1)} \int_{\fT^1} f(y) d \nu(y),
\end{equation}
for $\nu$ a.a. points $v \in \fT^1$ and every $\nu$-measurable and integrable $f: \fT^1 \mapsto \fR$.
%\item[2)] The measure $\nu$ is self-similar, in the sense that there exist an interval $I$, centered on $\pi/2$, and a number $0<t <1 $, such that
%$$\nu  = t \ \nu \circ T_{|I|} \quad {\rm on} \quad I,$$
%where $$T_{|I|}(x)={1 \over |I|} \left( x-{\pi \over 2} \right)+{\pi \over 2}.$$ 
%\end{itemize}
\end{conjecture}

To support the conjecture we have computed the left and right hand sides of the equation $(\ref{ergodic})$ for several monomials in $\cos$ and $\sin$. The approximation of the  density of the measure has been computed by dividing the circles into a collection of $N$  subintervals $[\theta_i,\theta_{i+1}]$ of equal length, and counting the relative number $K_i^n$ of angles $\alpha_\omega^{F^*}$, $\omega \in \{0,1\}^{12}$, in $\cO_n$ that lie in those intervals. 

We have collected the data in Tables $\ref{ErgodicTable1}-\ref{ErgodicTable4}$. The number reported in the tables is $|2 (L-R)/(L+R)|$, where
$$L={1 \over M} \sum_{k=0}^{M-1} f \left( \cA_k((0,0)) v \right), \quad R=\sum_{i=0}^{N-1}f\left({\theta_i+\theta_{i+1} \over 2}\right) K_i^n.$$

Here, we have used that the point $(0,0)$ is in the stable Cantor set (see \cite{GJ2} for details). The quantity $L$ has been computed for two different $v$'s.

\begin{table}[h]
\begin{center}
\begin{tabular}{|r|l|l|r|}
\hline
$\cO_n$& $N=1000$  & $N=5000$ & $N=15000$ \\ \hline
$n=12$ &  $2.5767 \times 10^{-3}$  & $3.0478 \times 10^{-3}$ & $3.0450 \times 10^{-3}$  \\
$n=14$ &  $2.0387 \times 10^{-3}$  & $2.0535 \times 10^{-3}$ & $2.0552 \times 10^{-3}$  \\
$n=16$ &  $2.5675 \times 10^{-4}$  & $9.8508 \times 10^{-5}$ & $1.9443 \times 10^{-4}$ \\
\hline
\end{tabular}
\end{center}
\caption{The relative difference of the left and right hand sides in $(\ref{ergodic})$ for the function $\sin^2$, $M=20000$, $v=(1,0)$.
}\label{ErgodicTable1}
\end{table}

\begin{table}[h]
\begin{center}
\begin{tabular}{|r|l|l|r|}
\hline
$\cO_n$& $N=1000$  & $N=5000$ & $N=15000$ \\ \hline
$n=12$ & $4.3071 \times 10^{-3}$ & $4.2824 \times 10^{-3}$ & $4.2809 \times 10^{-3}$ \\
$n=14$ & $2.9367 \times 10^{-3}$ & $2.9642 \times 10^{-3}$ & $2.9731 \times 10^{-3}$ \\
$n=16$ & $2.4732 \times 10^{-4}$ & $2.1161 \times 10^{-4}$ & $2.0896 \times 10^{-4}$ \\
\hline
\end{tabular}
\end{center}
\caption{The relative difference of the left and right hand sides in $(\ref{ergodic})$ for the function $\cos^4$, $M=20000$, $v=(1,0)$.
}\label{ErgodicTable2}
\end{table}

\begin{table}[h]
\begin{center}
\begin{tabular}{|r|l|l|r|}
\hline
$\cO_n$& $N=1000$  & $N=5000$ & $N=15000$ \\ \hline
$n=12$ &  $2.6290 \times 10^{-3}$  & $3.1001 \times 10^{-3}$ & $3.0973 \times 10^{-3}$  \\
$n=14$ &  $1.6695 \times 10^{-3}$  & $1.6836 \times 10^{-3}$ & $1.6854 \times 10^{-3}$  \\
$n=16$ &  $1.1324 \times 10^{-4}$  & $4.6830 \times 10^{-4}$ & $5.6424 \times 10^{-4}$ \\
\hline
\end{tabular}
\end{center}
\caption{The relative difference of the left and right hand sides in $(\ref{ergodic})$ for the function $\sin^2$, $M=20000$, $v=(0,1)$.
}\label{ErgodicTable3}
\end{table}

\begin{table}[h]
\begin{center}
\begin{tabular}{|r|l|l|r|}
\hline
$\cO_n$& $N=1000$  & $N=5000$ & $N=15000$ \\ \hline
$n=12
$ & $4.4239 \times 10^{-3}$ & $4.3992 \times 10^{-3}$ & $4.3978 \times 10^{-3}$ \\
$n=14$ & $2.0916 \times 10^{-3}$ & $2.1189 \times 10^{-3}$ & $2.1278 \times 10^{-3}$ \\
$n=16$ & $1.0927 \times 10^{-3}$ & $1.0568 \times 10^{-3}$ & $1.0542 \times 10^{-3}$ \\
\hline
\end{tabular}
\end{center}
\caption{The relative difference of the left and right hand sides in $(\ref{ergodic})$ for the function $\cos^4$, $M=20000$, $v=(0,1)$.
}\label{ErgodicTable4}
\end{table}

\section{Reality of the renormalization spectrum}

It has been a long standing conjecture, supported by numerical evidence, that the renormalization spectrum of the Feigenbaum period-doubling operator for unimodal maps is real. We observe the same phenomenon with the period-doubling renormalization spectrum for area-preserving maps

We recall that our Theorem $\ref{Main_Theorem_1}$ has been proved for the operator $\cR_{EKW}$ acting on the space $\cA_s(1.75)$.

%Recall the Definition  $\ref{B_space}$ of the Banach subspace $\cA_s(\rho)$ of $\cA(\rho)$. 
We will choose a new basis $\{\psi_{i,j}\}$ in $\cA_s(\rho)$. Given $s \in \cA_s(\rho)$ we write its Taylor expansion in the form
$$s(x,y)=\sum_{(i,j)\in I} s_{i,j} \psi_{i,j}(x,y),$$ 
where $\psi_{i,j} \in \cA_s(\rho)$:
\begin{eqnarray}
\nonumber \tilde{\psi}_{i,j}(x,y)&=&x^{i+1} y^j,  \quad i=-1, \quad j \ge 0, \\
\nonumber \tilde{\psi}_{i,j}(x,y)&=& x^{i+1} y^j +{i+1 \over {j+1}}  x^{j+1}  y^i,  \quad i >-1, \quad j \ge i,\\
\nonumber \psi_{i,j}&=&{\tilde{\psi}_{i,j} \over \|\tilde{\psi}_{i,j} \|_\rho}, \quad i \ge -1, \quad j \ge \max\{0,i\},
\end{eqnarray}
and the index set $I$ of these basis vectors is defined as
$$I=\{(i,j) \in \field{Z}^2: \quad i \ge -1, \quad j \ge \max\{0,i\}  \}.$$

%Denote $\tilde{\cA}_s(\rho)$ the set of all sequences 
%$$\tilde{s}=\left\{s_{i,j}: s_{i,j} \in \field{C}, \sum_{(i,j) \in I} |s_{i,j}|  < \infty \right\}.$$
%Equipped with the $l_1$-norm 
%\begin{equation}\label{l_1}
%|s|_1=\sum_{(i,j) \in I} |s_{i,j}|,
%\end{equation}
%$\tilde{\cA}_s(\rho)$ is a Banach space, which is isomorphic to $\cA_s(\rho)$. Clearly, the isomorphism $J: \cA_s(\rho) \mapsto \tilde{\cA_s}(\rho)$ is an isometry:
%$$\| \cdot \|_\rho=| \cdot |_1.$$

%We divide the set $I$  in three disjoint  parts: 
%\begin{eqnarray}
%\nonumber I_1 &=& \{(i,j) \in I: i+j<N \},\\
%\nonumber  I_2 &=& \{(i,j) \in I:  N  \le i+j < M \},\\
%\nonumber I_3 &=& \{(i,j) \in I:  i+j \ge M \},
%\end{eqnarray} 
%with
%$$N=22, \quad M=60.$$

We will denote the cardinality of the set
$$I_N=\{(i,j) \in I: i+j<N \}$$
as $D(N)$. Assign a single index to vectors $\psi_{i,j}$, $(i,j) \in I_N$, as follows:
\begin{eqnarray}
\nonumber &&k(-1,0)=1,\quad k(-1,1)=2, \quad  \ldots ,\quad  k(-1,N)=N+1, \quad k(0,0)=N+2, \\
\nonumber &&  k(0,1)=N+3,\quad \ldots, \quad k\left(\left[{N-1\over 2} \right], N-1-\left[{N-1\over 2} \right]\right)=D(N).
\end{eqnarray}
This correspondence $(i,j)  \mapsto k$ is one-to-one, we will, therefore, also use the notation $(i(k),j(k))$. 

%First, we compute the action of the renormalization derivative on each $\psi_{k(i,j)}$, $(i,j)\in I_1$. 

For any $s \in \cA_s(\rho)$, we define the following projections on the subspaces of the linear subspace $E_{D(N)}$ spanned by $\{\psi_k\}_{k=1}^{D(N)}$.
$$\Pi_k s=s_{i(k),j(k)} \psi_k, \quad \Pi_{E_{D(N)}} s = \sum_{m \le {D(N)}}  \Pi_m s.$$
%\begin{eqnarray}
%\nonumber \Pi_{D(N)}  \psi&=&\left( \psi \right)_{i(D(N))+1,j(D(N))} \|\tilde{\psi}_{D(N)} \|_\rho \psi_{D(N)},\\
%\nonumber \Pi_{D(N)-l}  \psi&=& \left( \psi -\sum_{m={D(N)}}^{D(N)-l+1}  \Pi_m \psi \right)_{i(D(N)-l)+1,j(D(N-l))} \!\!\!\!\!\!\!\!\!\!\!\!\!\!\!\!\!\!\!\!\!\!\!\!\!\!\!\!\!\!\!\!\!\!\!\! \|\tilde{\psi}_{D(N)-l} \|_\rho  \psi_{D(N)-l}, \quad 1 \le  l \le D(N)-1,\\
%\nonumber  \Pi_{E_{D(N)}} \psi &=& \sum_{m \le {D(N)}}  \Pi_m \psi.
%\end{eqnarray}
%It is easily verified that $\Pi_m$ are indeed projections, $\Pi_m^2=\Pi_m$.

We have found an approximation of a matrix representation of the finite-dimensional linear operator 
$$\Pi_{E_{D(N)}}  D \cR_{EKW}[s^*]  \Pi_{E_{D(N)}}$$
as 
$$ D_{n,m}= \Pi_m   D \cR_{EKW}[s^*]  \psi_n ,$$
and have computed the eigenvalues of the matrix. Table $\ref{eigenvalues_40}$ gives the first $100$ eigenvalues of the operator $\Pi_{E_{D(N)}}  D \cR_{EKW}[s^*]  \Pi_{E_{D(N)}}$, $N=20$.

We would like to note that the complexity of some of the higher-order eigenvalues is likely a numerical artifact. To demonstrate this, we list eigenvalues  $33$ and $34$ of the operator $\Pi_{E_{D(N)}}  D \cR_{EKW}[s^*]  \Pi_{E_{D(N)}}$ for several values of $N$ in Table $\ref{comparison1}$, and eigenvalues $62$ and $63$ in Table $\ref{comparison2}$. One can see that as one increases the dimension of the projection, the pair of complex conjugate eigenvalues turns real after some $N$, and stays that way.

\begin{center}

\begin{small}

\begin{table}[h]
\begin{center}
\begin{tabular}{|l|l|l|}
\hline
$N$ & \hspace{3.5cm}  $33$ & \hspace{3.5cm}$  34$\\
\hline
$10$ &  $7.710128558752339  \times 10^{-7}-1.061214359986916 \times 10^{-7} i$  &   $7.710128558752339  \times 10^{-7}+1.061214359986916 \times 10^{-7} i$ \\
\hline
$15$ & $9.116707201322168 \times 10^{-7}+0.0 i$ & $8.866211124768190 \times 10^{-7} +0.0 i$ \\
\hline
$20$ & $9.116282429256495 \times 10^{-7}+0.0 i$ &  $8.873952087973409 \times 10^{-7} +0.0 i$ \\
\hline
\end{tabular}
\end{center}
\caption{A comparison of the eigenvalues $33$ and $34$ of the operator $\Pi_{E_{D(N)}}  D \cR_{EKW}[s^*]  \Pi_{E_{D(N)}}$ for several values of $N$.}\label{comparison1}
\end{table}

\begin{table}[h]
\begin{center}
\begin{tabular}{|l|l|l|}
\hline
$N$ & \hspace{4cm} $62$ & \hspace{4cm}  $63$\\
\hline
$15$ &  $-6.705013355269772 \times 10^{-10}-1.625832021663626 \times 10^{-10} i$  &  $-6.705013355269772 \times 10^{-10}+1.625832021663626 \times 10^{-10} i$ \\
\hline
$20$ & $-8.704206290836820 \times 10^{-10}+0.0 i$ & $-8.227051108116014 \times 10^{-10} +0.0 i $ \\
\hline
$25$ & $-8.704155784717757 \times 10^{-10}+0.0 i$ & $-8.247546278128800 \times 10^{-10} +0.0 i $ \\
\hline
\end{tabular}
\end{center}
\caption{A comparison of the eigenvalues $62$ and $63$ of the operator $\Pi_{E_{D(N)}}  D \cR_{EKW}[s^*]  \Pi_{E_{D(N)}}$ for several values of $N$.}\label{comparison2}
\end{table}
\end{small}
\end{center}

\begin{small}
\begin{table}
\begin{minipage}{0.37\linewidth}

\begin{tabular}{|l|l|l|}
\hline
$ 1$  &  $ \phantom{-}8.72109720060341027+ 0.0 i $ &$N$ \\
$ 2$  &  $-4.01807670479890925+ 0.0 i $ & $Y$\\
$ 3$  &  $-2.48875288718523027 \times 10^{-1} +0.0 i $ & $Y$ \\
$ 4$  &  $-1.16629420927308277 \times 10^{-1}+ 0.0 i $  &$N$ \\
$ 5$  &  $ \phantom{-}7.29842918134375602 \times 10^{-2} + 0.0 i $ &$N$ \\
$ 6$  &  $ \phantom{-}6.19389093347281188 \times 10^{-2} + 0.0 i $ &$Y$ \\
$ 7$  &  $-1.54150639435907658 \times 10^{-2} + 0.0 i $ & $Y$ \\
$ 8$  &  $-1.50053025013312190 \times 10^{-2} + 0.0 i $ &$Y$ \\
$ 9$  &  $-5.24316259408011873 \times 10^{-3} + 0.0 i $ &$N$ \\
$ 10$  &  $ \phantom{-}3.83642848957628521 \times 10^{-3} + 0.0 i $ & $Y$ \\
$ 11$  &  $ \phantom{-}3.73444899232657950 \times 10^{-3} + 0.0 i $ &$Y$ \\
$ 12$  &  $ \phantom{-}1.67255839935358406 \times 10^{-3} + 0.0 i $ &$N$\\
$ 13$  &  $ \phantom{-}1.95630639839669629 \times 10^{-3} + 0.0 i $ &$N$\\
$ 14$  &  $-9.54792247990135266 \times 10^{-4} + 0.0 i $ &$Y$\\
$ 15$  &  $-9.29412071169940800 \times 10^{-4} + 0.0 i $ &$Y$\\
$ 16$  &  $-7.74683770828318898 \times 10^{-4} + 0.0 i $ &$N$\\
$ 17$  &  $ \phantom{-}2.37624196383790100 \times 10^{-4} + 0.0 i $&$Y$\\
$ 18$  &  $ \phantom{-}2.31307697553682750 \times 10^{-4} + 0.0 i $&$Y$\\
$ 19$  &  $ \phantom{-}2.05791964161655016 \times 10^{-4} + 0.0 i $&$N$\\
$ 20$  &  $ \phantom{-}4.61021338856759535 \times 10^{-5} + 0.0 i $&$N$\\
$ 21$  &  $-5.91387904855386636 \times 10^{-5} + 0.0 i $&$Y$\\
$ 22$  &  $-5.75667700266775337 \times 10^{-5} + 0.0 i $&$Y$\\
$ 23$  &  $-5.60365367884625166 \times 10^{-5} + 0.0 i $&$Y$\\
$ 24$  &  $-3.87086314650991711 \times 10^{-5} + 0.0 i $&$N$\\
$ 25$  &  $ \phantom{-}1.47181835587117418 \times 10^{-5} + 0.0 i $&$Y$\\
$ 26$  &  $ \phantom{-}1.43269464891488572 \times 10^{-5} + 0.0 i $&$Y$\\
$ 27$  &  $ \phantom{-}1.39461092667464458 \times 10^{-5} + 0.0 i $&$Y$\\
$ 28$  &  $ \phantom{-}3.37032628919152810 \times 10^{-6} + 0.0 i $&$N$\\
$ 29$  &  $-3.66299218723578256 \times 10^{-6} + 0.0 i $&$Y$\\
$ 30$  &  $-3.56562288094160436 \times 10^{-6} + 0.0 i $&$Y$\\
$ 31$  &  $-3.47084200041532543 \times 10^{-6} + 0.0 i $&$Y$\\
$ 32$  &  $-2.68212102338467744 \times 10^{-6} + 0.0 i $&$N$\\
$ 33$  &  $ \phantom{-}9.11628242925649478 \times 10^{-7} + 0.0 i $&$Y$\\
$ 34$  &  $ \phantom{-}8.87395208797340851 \times 10^{-7} + 0.0 i $&$Y$\\
$ 35$  &  $ \phantom{-}8.63807486025034232 \times 10^{-7} + 0.0 i $&$Y$\\
$ 36$  &  $ \phantom{-}8.43743910455452052 \times 10^{-7} + 0.0 i $&$N$\\
$ 37$  &  $-6.13340779558256008 \times 10^{-7} + 0.0 i $&$N$\\
$ 38$  &  $ \phantom{-}3.64333304618918514 \times 10^{-7} + 0.0 i $&$N$\\
$ 39$  &  $-2.26881741829139139 \times 10^{-7} + 0.0 i $&$Y$\\
$ 40$  &  $-2.20850766682014296 \times 10^{-7} + 0.0 i $&$Y$\\
$ 41$  &  $-2.14980452161630267 \times 10^{-7} + 0.0 i $&$Y$\\
$ 42$  &  $-2.09265346840381688 \times 10^{-7} + 0.0 i $&$Y$\\
$ 43$  &  $ \phantom{-}5.64652471456064890 \times 10^{-8} + 0.0 i $&$Y$\\
$ 44$  &  $ \phantom{-}5.49645085465527861 \times 10^{-8} + 0.0 i $&$Y$\\
$ 45$  &  $ \phantom{-}5.35030991885504657 \times 10^{-8} + 0.0 i $&$Y$\\
$ 46$  &  $ \phantom{-}5.20806572338405750 \times 10^{-8} + 0.0 i $&$Y$\\
$ 47$  &  $ \phantom{-}4.06666476717650202 \times 10^{-8} + 0.0 i $&$N$\\
$ 48$  &  $-2.89661341440091808 \times 10^{-8} + 0.0 i $&$N$\\
$ 49$  &  $-1.44725130238278229 \times 10^{-8} + 0.0 i $&$N$\\
$ 50$  &  $-1.40530550882139069 \times 10^{-8} + 0.0 i $&$Y$\\
\hline
\end{tabular}
\end{minipage}
\hspace{0.0001cm}
\begin{minipage}{0.61\linewidth}

\begin{tabular}{|l|l|l|}
\hline
$ 51$  &  $-1.36773104970236911 \times 10^{-8} + 0.0 i $&$Y$\\
$ 52$  &  $-1.29617185901792936 \times 10^{-8} + 0.0 i $ &$Y$\\
$ 53$  &  $-1.33177097864396282 \times 10^{-8} + 0.0 i $&$Y$\\
$ 54$  &  $ \phantom{-}5.83632841811164741 \times 10^{-9} + 0.0 i $&$N$\\
$ 55$  &  $ \phantom{-}3.49734698280496765 \times 10^{-9} + 0.0 i $&$Y$\\
$ 56$  &  $ \phantom{-}3.40553232521175106 \times 10^{-9} + 0.0 i $&$N$\\
$ 57$  &  $ \phantom{-}3.30875856445169904 \times 10^{-9} + 0.0 i $&$N$\\
$ 58$  &  $ \phantom{-}3.23415716696474633 \times 10^{-9} + 0.0 i $&$N$\\
$ 59$  &  $ \phantom{-}3.20031850876952441 \times 10^{-9} + 0.0 i $&$N$\\
$ 60$  &  $-1.21497298394740393 \times 10^{-9} + 0.0 i $&$N$\\
$ 61$  &  $-8.47423932266715023 \times 10^{-10} + 0.0 i $&$Y$\\
$ 62$  &  $-8.70420629083682033 \times 10^{-10} + 0.0 i $&$Y$\\
$ 63$  &  $-8.22705110811601487 \times 10^{-10} + 0.0 i $&$N$\\
$ 64$  &  $-8.07179147279961451 \times 10^{-10} + 0.0 i $&$N$\\
$ 65$  &  $-7.79143693124703818 \times 10^{-10} + 0.0 i $&$N$\\
$ 66$  &  $ \phantom{-}4.63752609505868076 \times 10^{-10} + 0.0 i $&$N$\\
$ 67$  &  $-3.22686335036190861 \times 10^{-10} + 0.0 i $&$N$\\
$ 68$  &  $ \phantom{-}2.77581537948968217 \times 10^{-10} + 0.0 i $&$N$\\
$ 69$  &  $ \phantom{-}2.16643085621470110 \times 10^{-10} + 0.0 i $&$Y$\\
$ 70$  &  $ \phantom{-}2.09856849656907084 \times 10^{-10} + 0.0 i $&$N$\\
$ 71$  &  $ \phantom{-}2.08109796932933800 \times 10^{-10} + 0.0 i $&$N$\\
$ 72$  &  $ \phantom{-}1.96066223947378410 \times 10^{-10} -4.43072486498362171 \times 10^{-12} i $&$N$\\
$ 73$  &  $ \phantom{-}1.96066223947378410 \times 10^{-10} + 4.43072486498362171 \times 10^{-12} i $&$N$\\
$ 74$  &  $-4.00724308724832731 \times 10^{-11} -1.29484585578379769 \times 10^{-11} i $&$N$\\
$ 75$  &  $-4.00724308724832731 \times 10^{-11} + 1.29484585578379769 \times 10^{-11} i $&$N$\\
$ 76$  &  $-5.39409517510675544 \times 10^{-11} + 0.0 i $&$Y$\\
$ 77$  &  $-5.22615941114597117 \times 10^{-11} -6.44631270712982928 \times 10^{-13} i $&$N$\\
$ 78$  &  $-5.22615941114597117 \times 10^{-11} + 6.44631270712982928 \times 10^{-13} i $&$N$\\
$ 79$  &  $-4.83301369903105249 \times 10^{-11} -1.61011502564644013 \times 10^{-12} i $&$N$\\
$ 80$  &  $-4.83301369903105249 \times 10^{-11} + 1.61011502564644013 \times 10^{-12} i $&$N$\\
$ 81$  &  $-1.21652373238263412 \times 10^{-11} + 0.0 i $&$N$\\
$ 82$  &  $ \phantom{-}1.34403330160617114 \times 10^{-11} + 0.0 i $&$Y$\\
$ 83$  &  $ \phantom{-}1.29886814444646485 \times 10^{-11} -4.49484074640346900 \times 10^{-13} i $&$N$\\
$ 84$  &  $ \phantom{-}1.29886814444646485 \times 10^{-11} + 4.49484074640346900 \times 10^{-13} i $&$N$\\
$ 85$  &  $ \phantom{-}1.30874465804588329 \times 10^{-11} + 0.0 i $&$N$\\
$ 86$  &  $ \phantom{-}1.22553839952396333 \times 10^{-11} + 0.0 i $&$N$\\
$ 87$  &  $ \phantom{-}1.14128486730559826 \times 10^{-11} -1.22675413039504159 \times 10^{-12} i $&$N$\\
$ 88$  &  $ \phantom{-}1.14128486730559826 \times 10^{-11} + 1.22675413039504159 \times 10^{-12} i $&$N$\\
$ 89$  &  $-4.64857574616255465 \times 10^{-13} -4.15935498398141420 \times 10^{-12} i $&$N$\\
$ 90$  &  $-4.64857574616255465 \times 10^{-13} + 4.15935498398141420 \times 10^{-12} i $&$N$\\
$ 91$  &  $-3.50139347550126033 \times 10^{-12} + 0.0 i $&$N$\\
$ 92$  &  $-3.12844287139405750 \times 10^{-12} -2.35599032104882566 \times 10^{-13} i $&$Y$\\
$ 93$  &  $-3.12844287139405750 \times 10^{-12} + 2.35599032104882566 \times 10^{-13} i $&$Y$\\
$ 94$  &  $-3.18763936883361579 \times 10^{-12} + 0.0 i $&$N$\\
$ 95$  &  $-2.60876535503609472 \times 10^{-12} + 0.0 i $&$N$\\
$ 96$  &  $ \phantom{-}1.00716118446860206 \times 10^{-12} + 0.0 i $&$N$\\
$ 97$  &  $-8.13328302601060092 \times 10^{-13} + 0.0 i $&$N$\\
$ 98$  &  $ \phantom{-}8.04488862233887989 \times 10^{-13} -1.00262940606313244 \times 10^{-13} i $&$N$\\
$ 99$  &  $ \phantom{-}8.04488862233887989 \times 10^{-13} + 1.00262940606313244 \times 10^{-13} i $&$N$\\
$ 100$  &  $ \phantom{-}5.08237476178472880 \times 10^{-13} + 0.0 i $&$Y$\\
%$ 100$  &  $ \phantom{-}4.84975108680808324 \times 10^{-13} -1.05606992205910460 \times 10^{-13} i $ &$N$\\
\hline
\end{tabular}
\end{minipage}

\bigskip

\caption{Some of the eigenvalues of $\Pi_{E_{D(N)}}  D \cR_{EKW}[s^*]  \Pi_{E_{D(N)}}$, $N=20$. The last column lists those eigenvalues that correspond to the products of powers of $\lambda_*$ and $\mu_*$.}\label{eigenvalues_40}
\end{table}

\end{small}

\section*{Acknowledgment}
The second author is funded by a postdoctoral fellowship from \textit{Vetenskapsr\aa det} (the Swedish Research Council).

\end{document}